\documentclass[12pt]{article}

\usepackage{amsmath,amssymb}
\usepackage{setspace,color,enumitem,url}

\begin{document}

\begin{center}
\textbf{\Large On valid descriptive inference from\\ non-probability sample}\\
\emph{Li-Chun Zhang}\footnote{Address for correspondence: S3RI/University of Southampton, Highfield SO17 1BJ, Southampton, UK. Email: L.Zhang@soton.ac.uk}
\end{center}

\emph{Abstract:} We examine the conditions under which descriptive inference can be based directly on the observed distribution in a non-probability sample, under both the super-population and quasi-randomisation modelling approaches. Review of existing estimation methods reveals that the traditional formulation of these conditions may be inadequate due to potential issues of under-coverage or heterogeneous mean beyond the assumed model. We formulate unifying conditions that are applicable to both type of modelling approaches. The difficulties of empirically validating the required conditions are discussed, as well as valid inference approaches using supplementary probability sampling. The key message is that probability sampling may still be necessary in some situations, in order to ensure the validity of descriptive inference, but it can be much less resource-demanding provided the presence of a big non-probability sample.  


\section{Introduction}

There is a resurgence of interest in the use of non-probability samples. See e.g. Baker et al. (2013) and Elliot and Valliant (2017) for two recent reviews. Such data may arise in situations where probability sampling is either infeasible or too costly. The observations may be obtained from the so-called big-data sources, such as payment transaction data via a specific platform, cellphone call data from a major provider of the service. These big-data non-probability samples can be much larger in size, compared to the more familiar non-probability samples collected from web panel surveys, quota sampling, etc.  

Following Rubin (1976) and Little (1982), Smith (1983) consider the so-called \emph{super-population (SP)} approach to inference from non-probability sample. Under this approach, a prediction model is constructed for the outcome variable of interest, often conditional on some chosen covariates. In particular, Smith (1983) observes an important distinction between \emph{analytic} and \emph{descriptive} inference. In analytic inference, the target is the prediction model parameters that are of a theoretical nature; such parameters can never be observed directly no matter how large the sample is. Whereas the targets of descriptive inference are statistics of a given finite population, such that in principle they can be directly observed provided a perfect census of the population.
 
Moreover, Smith (1983) focuses on \emph{validity conditions}, under which the non-probability sample observation mechanism can be ignored, in the sense that inference can be based on the \emph{observed} distributions directly, such as the conditional distribution of the outcome variable given the covariates in the sample. The two key validity conditions under the SP approach can be roughly stated as follows: (i) the prediction model is correctly specified for the population units, (ii) the non-probability sample selection mechanism is non-informative, in the sense that the relevant distribution under the population model can be observed in the non-probability sample directly. Similar validity conditions for the SP approach apply in other situations, such as purposive sampling (Royall, 1970), missing data problems (Rubin, 1976). 

In this paper we concentrate on descriptive inference methods that depend on validity conditions in the sense of Smith (1983). Of course, inference is also possible without such validity conditions. For instance, not missing-at-random models (Rubin, 1976) can be used to deal with informative missing data, where the unobserved full-sample outcome distribution is not the same as that among the respondent subsample. Or, the sample likelihood of Pfeffermann et al. (1998) can be applied to survey data under informative sampling, where the distribution that holds in the population cannot be directly observed in the sample. See also Pfeffermann (2017) for several other situations where this approach may be relevant. However, in this paper we do not consider such methods, which explicitly address informative observation mechanisms of sample selection or measurement. 

As reviewed by Elliot and Valliant (2017), there exists another \emph{quasi-randomisation (QR)} approach to non-probability samples. Under the QR approach, one hypothesises a randomisation model of the non-porbability sample inclusion indicator, but treats the outcomes of interest as unknown constants in the population. Though it is clearly inspired by the randomisation approach based on probability sampling, the QR approach is also a model-based approach, based on a model of the sample inclusion indicator instead of a prediction model of the outcome variable under the SP approach. A key motivation is that the correct inclusion probability can be used for any outcome of interest, just like when it is known under probability sampling, whereas the SP approach by nature must be specified differently for different outcome variables. In the context of survey sampling, the QR approach was introduced to deal with nonresponse, where response to survey is modelled as the second phase of selection, in addition to the first phase of sample selection according to a probability sampling design (Oh and Scheuren, 1983). 

According to Elliot and Valliant (2017), two key validity conditions are required for the QR approach. (I) The non-probability sample does have a probability sampling mechanism, even though it is unknown. In particular, one assumes that this hypothesised sample inclusion probability is strictly positive for all the population units, so that the only difference to probability sampling is that the inclusion probability is unknown. (II) There exist a set of covariates that ``fully govern the sampling mechanism''. In other words, the sample inclusion probability is a function of these covariates.  
  
Thus, there are two model-based approaches to inference from non-probability sample. Under the SP approach, one models the outcome variable conditional on the realised sample inclusion indicators; whereas under the QR approach, one models the sample inclusion indicators, but treats the outcomes as unknown constants. Although one may envisage the outcomes as the realised values of random variables, a fully specified model of the outcome variable will not be required under the QR approach, provided suitable validity conditions. Similarly, although one acknowledges that the sample selection mechanism may be critical to the SP approach, a fully specified model of the inclusion indicator will not be required under the SP approach, provided suitable validity conditions. 
  
It is possible to construct estimators that combine both the models of outcome and sample inclusion indicator, in a manner such that the estimator is consistent as long as one of the two models hold, \emph{provided} the same covariates are used in both models. Over the recent years, it is becoming common to refer to this estimation approach as ``doubly robust''.  Still, in reality, how likely is it for both the sample inclusion mechanism and the outcome generation mechanism to be fully explained by exactly the same covariates? How likely is this to be the case moving between different outcomes? Notice that the traditional generalised regression estimator in survey sampling is doubly-robust in the same sense. It is a fact that in the debate between model-based vs. design-based inference from probability sampling, either side questioned the ``robustness'' of the other approach. 

The rest of the paper is organised as follows. In Section \ref{review} we review the estimation methods for non-probability sample which do require validity conditions. Although these have been roughly stated above, a closer examination under both modelling perspectives reveals nuances across the different estimators. Moreover, we shall highlight the potential challenges of under-coverage and heterogeneous means beyond the assumed model. The traditional formulation of validity conditions is inadequate in both regards. We outline a set of unified validity conditions in Section \ref{condition}, which are formulated non-parametrically and encompasses both the modelling approaches. Post-stratification and calibration estimators are considered in light of these conditions. However, as will be discussed, a key difficulty in practice is that the validity conditions may be impossible to verify empirically based only on the data used for the estimation. Finally, we outline shortly in Section \ref{combine} two approaches given a supplementary probability sampling of the outcome of interest. The key message is that probability sampling may still be necessary in some situations, in order to ensure the validity of descriptive inference, but it can be much less resource-demanding provided the presence of a big non-probability sample. In fact, the bigger the non-probability sample, the better it is.

\section{Review of existing approaches} \label{review}

Denote by $U$ the population of known size $N$. Let each population unit be associated with an  \emph{outcome} of interest, denoted by $y_i$, for $i\in U$. Denote by $B$ the \emph{observed} nonprobability sample of size $n_B$. A common assumption to all the estimators we discuss below is that $B$ does not contain any out-of-scope units, such that $B\subset U$. Let $\delta_i = 1$ if $i\in B$, and 0 if $i\in U\setminus B$. Let $y_i$ be observed for all the units in $B$, and let $y_B = \{ y_i; i\in B\}$. To fix the idea, let 
\[
Y = \sum_{i\in U} y_i
\] 
be the population total that is the target of descriptive inference. 
Let $x_B = \{ x_i; i\in B\}$ in cases where any relevant covariates $x_i$ are available in the sample $B$. Let $X = \sum_{i\in U} x_i$ be the population totals and let $\bar{X} = X/N$. Provided $x_B$ is available, one can have two situations depending on whether $(X, \bar{X})$ are known or not. In the case they are unknown, it may still be possible that there exists a second \emph{probability sample} $S$, for $S\subset U$, in which $x_i$ is observed, so that $(X, \bar{X})$ can be estimated based on the sample $S$. 


\subsection{$B$-sample expansion estimator}

Consider first the most basic situation where only $y_B$ is observed, and no relevant covariates are available at all. The $B$-sample expansion estimator of $Y$ is given by
\begin{equation}
\widehat{Y} = N \bar{y}_B \label{mean}
\end{equation}
where $\bar{y}_B = \sum_{i\in B} y_i/n_B$ is the $B$-sample mean. 

Under the SP approach, let 
\[
\mu_i = E(y_i | \delta_i, i\in U)
\]
be the conditional expectation of $y_i$ given $\delta_i$, for any $i\in U$, where both $\delta_i$ and $y_i$ are treated as random variables. Provided the conditional expectation is the same as the unconditional expectation, for any $i\in U$, denoted by
\begin{equation}
\mu_i = \mu \label{SP}
\end{equation}
we have
\[
E(\bar{y}_B - Y/N | B) = \sum_{i\in B} \mu/n_B - \mu = 0
\] 
such that the $B$-sample expansion estimator is prediction unbiased for $Y$. We shall refer to \eqref{SP} as the \emph{SP assumption}, which is a validity condition for the $B$-sample expansion estimator under the SP approach.

Under the QR approach, where $y_i$ is treated as a fixed constant, let
\[
p_i = \mbox{Pr}(\delta_i = 1; y_i, i\in U)
\]
be the \emph{inclusion probability} of any population unit that is associated with the value $y_i$. The notational difference between ``;'' and ``$|$'' is introduced because, strictly speaking, $p_i$ is not a conditional probability now that $y_i$ is not conceived as the realised value of a random variable under the QR approach. Now, provided the inclusion probability is the same for any $i\in U$, denoted by
\begin{equation}
p_i = p \label{QR}
\end{equation}
we have $\widetilde{Y} = \sum_{i\in B} y_i/p$ is unbiased for $Y$, since
\[
E(\sum_{i\in B} y_i/p) = \sum_{i\in U} E(\delta_i; y_i, i\in U) y_i/p = \sum_{i\in U} p y_i/p = Y
\]
In reality, $p$ is unknown. Under \eqref{QR}, it is natural to estimate it by $\hat{p} = n_B/N$, which yields \eqref{mean} as the resulting plug-in estimator. It follows that the \emph{QR assumption} \eqref{QR} is the key validity condition, which ensures that the $B$-sample expansion estimator is consistent for $Y$, as $N\rightarrow \infty$ and $n_B/N = O_p(1)$ asympotically. 

In summary, the $B$-sample expansion estimator \eqref{mean} can be motivated under both the SP and QR approaches, provided validity conditions \eqref{SP} and \eqref{QR}, respectively.

\subsection{$B$-sample calibration estimator}

Suppose relevant covariates $x_B$ are available in the sample $B$. The population totals $X$ may be either known or unknown. In the latter case, suppose they can be estimated from a second probability sample $S$. The $B$-sample calibration estimator of $Y$ is given by
\begin{align}
& \widehat{Y} = \sum_{i\in B} w_i y_i \qquad\text{where} \quad
\begin{cases}  \sum_{i\in B} w_i x_i = X & \text{if known } X \\
\sum_{i\in B} w_i x_i = \widehat{X}(S) & \text{if unknown } X \end{cases} \label{cal} 
\end{align}
where $\widehat{X}(S)$ is some consistent $S$-sample estimator, as the $S$-sample size increases, and the weights $w_B = \{ w_i; i\in B\}$ are calibrated in a way depending on the availability of $X$. 
 
To actually compute the estimator \eqref{cal}, one needs to choose a set of initial weights, denoted by $a_B = \{ a_i; i\in B\}$. In the case of 
\begin{equation}
a_i = 1/p_i \label{QRcal}
\end{equation}
where $p_i$ is the true $B$-sample inclusion probability, for $p_i > 0$, the calibration estimator is asymptotically consistent, as $N\rightarrow \infty$ and $n_B/N = O_p(1)$, provided mild regularity conditions in addition. However, insofar as one cannot manage to set the initial weights \eqref{QRcal}, the calibration estimator is unmotivated from the QR perspective.  

Next, under the SP approach, suppose the SP$_x$ assumption given by
\begin{equation}
E(y_i | x_i, i\in U) = \mu(x_i) = x_i^{\top} \beta \label{SPx}
\end{equation}
which relates the conditional expectation of $y_i$ linearly to the given $x_i$, and
\begin{equation}
E(y_i | x_i, i\in U) = E(y_i | \delta_i, x_i, i\in U) = E(y_i | x_i, i\in B) \label{niy}
\end{equation}
by which the $B$-sample selection is \emph{non-informative} given $x_i$. We have then
\[
E(\sum_{i\in B} w_i y_i - Y | x_U) = E(\sum_{i\in B} w_i x_i^{\top} \beta) - X^{\top} \beta = 0
\]
provided $\sum_{i\in B} w_i x_i = X$, regardless of the initial weights $a_B$. Otherwise, this expectation would tend to 0, provided is $\widehat{X}(S)$ is an asymptotically unbiased estimator of $X$, under some suitable asymptotic setting. It follows that the assumptions \eqref{SPx} and \eqref{niy} are the key validity conditions for the $B$-sample calibration estimator under the SP approach.

The estimator \eqref{cal} becomes the $B$-sample post-stratification estimator in the special case where $x_i$ is the post-stratum dummy index. For the QR approach, one can set $a_i$ to be the inverse post-stratum $B$-sample fraction, which is equivalent to introducing the QR assumption \eqref{QR} in each post-stratum separately. This \emph{QR$_x$ assumption} provides then a validity condition for the $B$-sample post-stratification estimator under the QR approach. For the SP approach, the two assumptions \eqref{SPx} and \eqref{niy} remain formally the same.

\subsection{$B$-sample inverse propensity weighting}

Suppose relevant covariates $x_B$ are available in the sample $B$. The $B$-sample inverse propensity weighting (IPW) estimator is constructed under the QR approach. Suppose
\begin{equation}
p_i = p(x_i; \eta) >0 \label{QRx}
\end{equation}
i.e. the $B$-sample inclusion probability is completely determined given $x_i$, in the strictly positive parametric form $p(x_i; \eta)$, which may as well be referred to as the \emph{QR$_x$ assumption}. Provided $x_i$ is known for all the units in the population, $\eta$ can be estimated, say, by a population estimating equation 
\[
\sum_{i\in U} H(\delta_i; \eta) = 0
\]
where $E[H(\delta_i; \eta)] = 0$. Otherwise, suppose $x_S$ is observed in a second probability sample $S$, one can use the pseudo population estimating equation
\[
\sum_{i\in S} d_i H(\delta_i; \eta) = 0
\] 
(Kim and Wang, 2018), where $d_i$ is the sampling weight, for $i\in S$, or some $S$-sampling design-consistent adjustment of it. To ensure that $H(\delta_i; \eta)$ is the same in both of these two estimating equations, i.e. whether $i\in S$ or just $i\in U$, one needs to assume that $S$-sampling from $U$ is non-informative for $\delta_i$, so that
\begin{equation}
\mbox{Pr}(\delta_i = 1 | x_i, i\in S) = \mbox{Pr}(\delta_i = 1 | x_i, i\in U) \label{nid}
\end{equation}
Notice that, provided non-informativeness \eqref{nid}, one can also simply use the unweighted $S$-sample estimating equation, which is given by
\[
\sum_{i\in S} H(\delta_i; \eta) = 0
\]
instead of the pseudo population estimating equation. Having obtained the parameter estimate $\hat{\eta}$, one obtains $\hat{p}_i = p(x_i; \hat{\eta})$ and the $B$-sample IPW estimator
\begin{equation}
\widehat{Y} = \sum_{i\in B} y_i/ \hat{p}_i \label{IPW}
\end{equation}
which is consistent for $Y$ under mild regularity conditions, provided $\hat{\eta}$ is consistent for $\eta$ under some suitable asymptotic setting. It follows that the QR$_x$ assumption \eqref{QRx} is its key validity condition, whereas the non-informativeness assumption \eqref{nid} is needed in addition when $x_i$ is only available in a probability sample $S$ instead of the population.

\subsection{Another $B$-sample IPW estimator}

Elliot and Valliant (2017) discuss another IPW estimator \eqref{IPW}, where $p_i$ is obtained with the help of a second so-called reference probability sample $S$, and is given by
\begin{equation}
p_i \propto \mbox{Pr}(S_i =1 | x_i, i\in U) \frac{\mbox{Pr}(\delta_i =1 |x_i, i\in B\cup S)}{\mbox{Pr}(S_i =1 |x_i, i\in B\cup S)}  \label{ref}
\end{equation}
where $S_i = 1$ if $i\in S$ and 0 if $i\in U\setminus S$, and to fix the idea one may suppose $S\cap B = \emptyset$. Firstly, the QR$_x$ assumption \eqref{QRx} is retained. The definition of $p_i$ by \eqref{ref} can then be motivated as follows: 
\begin{align*}
\frac{\mbox{Pr}(\delta_i =1 | x_i, i\in U)}{\mbox{Pr}(S_i =1 | x_i, i\in U)} & \propto 
\frac{\mbox{Pr}(x_i | \delta_i =1, i\in U)}{\mbox{Pr}(x_i | S_i =1, i\in U)} 
\qquad\quad \Big[ \text{prop. to } \frac{\mbox{Pr}(\delta_i =1 | i\in U)}{\mbox{Pr}(S_i =1 | i\in U)} \Big] \\
& \propto \frac{\mbox{Pr}(x_i | \delta_i =1, i\in B\cup S)}{\mbox{Pr}(x_i | S_i =1, i\in B\cup S)} \\
& \propto \frac{\mbox{Pr}(\delta_i =1 | x_i, i\in B\cup S)}{\mbox{Pr}(S_i =1 | x_i, i\in B\cup S)}
\quad \Big[ \text{prop. to } \frac{\mbox{Pr}(\delta_i =1 | i\in B\cup S)}{\mbox{Pr}(S_i =1 | i\in B\cup S)} \Big]
\end{align*}
provided the $S$-sample inclusion probability is also fully determined by the same $x_i$ in the sense of \eqref{QRx}. In other words, the validity condition for the IPW estimator \eqref{IPW} based on \eqref{ref} is that the QR$_x$ assumption \eqref{QRx} holds for both the $B$-sample and the $S$-sample, \textit{given the same $x_i$}.

We make two observations. Firstly, despite the superficial resemblance to the propensity scoring method of Rosenbaum and Rubin (1983), the above argument for $p_i$ is not the same. As Rosenbaum and Rubin (1983) state clearly before their first enumerated equation, ``In this paper, the $N$ units in the study are viewed as a simple random sample from some population'', where $N$ is the size of the combined sample of treatment and non-treatment. The analogy to this combined sample here is $B\cup S$. However, it is generally untenable that $B\cup S$ can be treated as a simple random sample from the population. Secondly, for any given probability sample $S$, it is possible to identify the variables that determine the designed inclusion probability, denoted by $\pi_i = \pi(z_i)$, for $i\in U$. There arises thus a question, ``what if $\pi(z_i)$ differs considerably from $p(x_i, \hat{\eta})$?'' Moreover, one may have more than one probability sample in which $x_i$ is observed. There arises then a question, ``which reference sample should one use?''

\subsection{Sample matching estimator}

Rivers (2007) applies the SP approach in situations where a second probability sample $S$ is available. Replace the linear SP$_x$ assumption \eqref{SPx} by the SP$_x$ assumption below:
\begin{equation}
\| E(y_i | x_i, i\in U) - E(y_j | x_j, j\in U)\| = O\big( \|x_i - x_j\| \big) \label{SPxC}
\end{equation} 
provided suitable choice of the metric $\| \cdot \|$, as $N\rightarrow \infty$. Moreover, retain the non-informativeness assumption \eqref{niy} for the $B$-sample, such that the SP$_x$ assumption \eqref{SPxC} holds in the $B$-sample, provided $n_B/N = O_p(1)$ as $N\rightarrow \infty$. Assume that the same $x_i$ is also observed in $S$. The sample matching (SM) estimator is given by
\begin{equation}
\widehat{Y} = \sum\limits_{i\in S} d_i \hat{y}_i \label{SM}
\end{equation}
where $\hat{y}_i = y_{k_i}$, for $k_i = \arg \min\limits_{j\in B} \|x_i -x_j\|$, i.e. $y_{k_i}$ is the nearest-neighbour (NN) imputation value from the $B$-sample for $i\in S$.  

Assume first exact matching, where $x_{k_i} = x_i$ for all $i\in S$ and $k_i \in B$. We have then
\begin{align*}
E\big[ \sum_{i\in S} d_i E(\hat{y}_i | x_i) \big] & = E\big[ \sum_{i\in S} d_i E(y_i | x_i, i\in B) \big] = E\big[ \sum_{i\in S} d_i E(y_i | x_i, i\in U) \big] \\
& = \sum_{i\in U} E(y_i | x_i, i\in U) = E(Y | x_U)
\end{align*}
such that the SM estimator \eqref{SM} is prediction unbiased for $Y$. Notice that in the case of $S =U$, the SM estimator is just an NN-imputation method, which is prediction unbiased provided exact matching for $S$. Moreover, whether $S=U$ or not, the SM estimator will be less efficient than the prediction-imputed SM estimator 
\[
\hat{Y} = \sum_{i\in S} d_i E\big( x_i ; \widehat{\beta}(B)\big)
\]
whenever a correct parametric specification of the conditional mean (via $\beta$) is possible. Next, it is not difficult to see that the consistency of the SM estimator \eqref{SM} can be established under mild conditions, if one assumes asymptotic exact matching instead, i.e. 
\begin{equation}
\| x_i - x_{k_i} \| \rightarrow 0 \text{ in probability,} \label{match}
\end{equation}
for any $i\in S$, as $N\rightarrow \infty$ and $n_B/N = O_p(1)$. It follows that the assumptions \eqref{niy}, \eqref{SPxC} and \eqref{match} are the key validity conditions for the consistency of the SM estimator \eqref{SM}.

We make two observations. Firstly, an attractive feature of the NN-imputation is that the imputed sample $S$ looks more realistic and natural than, say, by the prediction imputation. However, unless the $S$-sampling is non-informative, the NN-imputed $S$-sample will not resemble the true $S$-sample that could have been observed, since
\[
E(\hat{y}_i | x_i, i\in S) = E(y_i | x_i, i\in U) \neq E(y_i | x_i, i\in S)
\]
where the inequality is the case unless $S$-sampling is non-informative in the sense of \eqref{niy}. Secondly, for any other covariate $z_i \neq x_i$, including when $z_i$ contains the $S$-sample design variables, we have 
\[
E(\hat{y}_i | z_i, x_i, i\in S) = E(y_i | x_i, i\in U) \neq E(y_i | z_i, x_i, i\in U) 
\] 
unless $y_i$ and $z_i$ are conditionally independent of each other given $x_i$, for $i\in U$. This is a general problem for statistical matching of variables associated with distinct units, i.e. $y_i$ associated with $x_i$ for some $i\in B$ and $z_i$ associated with the same value $x_i$ but for some different unit in $S$. The following example illustrates both remarks above. 

\paragraph{\em Example:} Let $y_i$ be independent of  $x_i \sim \mbox{Unif}(0,1)$, for any $i\in U$. Then, the SP$_x$ assumption \eqref{SPxC} holds trivially, as long as the marginal expectaion $E(y_i)$ exists. Next, suppose simple random sample $B$, so that the non-informative assumption \eqref{niy} holds, and $E(\hat{y}_i | x_i, i\in S) = E(y_i | i\in U)$ regardless of the exact matching assumption.
Suppose stratified simple random $S$-sampling with two strata of different sampling fractions, so that the $S$-sample inclusion probability is not a constant. Then, the $S$-sampling is informative (given $x_i$) as long as the population stratum means are different, since 
\[
E(\bar{y}_S | x_S, S) = E(\bar{y}_S | S) \neq E(\bar{Y}) = E(\bar{Y} | x_U)
\]
where $\bar{y}_S$ is the true $S$-sample mean that is unknown, since $y_i$ is not observed in $S$. It follows that the NN-imputed $S$-sample $\{ \hat{y}_i; i\in S\}$ would look like a sample generated by simple random sampling, rather than the actual stratified sampling. Moreover, the SM-estimator of stratum means, corresponding to say $z_i  = 1, 2$, respectively, will be biased for the population stratum means. $\square$

\subsection{Summary and discussion}

All the estimators from non-probability sample observations reviewed above are model-based, whether the modelling is carried out under the SP or QR approach. Two features regarding the model covariate $x_i$, for $i\in U$, are worth recapitulating:
\begin{itemize}[leftmargin=6mm]
\item compared to the situation with known $x_U$, making use of an additional probability sample $x_S$ entails a loss of efficiency, as can be expected; 

\item the availability of an additional probability sample without the outcome variable is not a principle advantage, since it does not simplify  the validity conditions compared to the situation where $x_U$ is known, but it does resolve the practical difficulty when $x_U$ is unavailable yet some functions of $x_U$ are needed for descriptive inference. 
\end{itemize}

As noted by Kim and Rao (2018), there is an important issue which does not appear to have received sufficient attention in the existing approaches, namely the potential under-coverage of the $B$-sample, when some of the units have in fact zero chance of being included in it. Under the SP approach, the existence of such units means that extrapolation of the estimated (conditional) distribution of $y_i$ in the $B$-sample to these units can only be based on subjective believes rather empirical evidence. The QR approach is equally affected, since inference built on the basis of randomisation would have fallen apart even when the probability $p_i$ were known for all the $B$-sample units, for which $p_i$ is strictly positive by definition, let alone when it is unknown and needs to be estimated. 

To address the potential under-coverage, Kim and Rao (2018) consider a two-phase SM estimator. Let the $S$-sample be partitioned into $S_1$ and $S_0$, such that $S_1 = \{ i; p_i > 0\}$ and $S_0 = \{ i; p_i = 0\}$. First, estimate this unobserved partition via the $B$-sample support: 
\[
\hat{S}_1 = \{ i; \min_{j\in B} \| x_i - x_j \| < \epsilon \}
\] 
Each $S$-sample unit that is unsupported in the $B$-sample $\epsilon$-neighbourhood is assigned to $\hat{S}_0$. Suppose this partition estimator is consistent in the following sense:
\[
| \hat{S}_1 \cup S_1 | - |\hat{S}_1 \cap S_1| \rightarrow 0 \text{ in probability,}
\]
as $N\rightarrow \infty$ and $\epsilon\rightarrow 0$. Next, the two-phase SM estimator is given as
\[
\widehat{Y} = \sum_{i\in \hat{S}_1} d_i w_{2i} \hat{y}_i 
\]
where $\sum_{i\in \hat{S}_1} d_i w_{2i} x_i = \sum_{i\in S} d_i x_i$. In other words, the under-coverage is dealt with by the calibration of the weights $w_{2i}$. This can be motivated, provided the conditional mean $E(y_i | x_i, p_i =0)$ can be linearly related to $x_i$, and the relationship is the same for the units with $p_i > 0$, i.e. the under-coverage is non-informative for the SP linear model. 
 
Now that inference from non-probability sample need to be model-based, there is the general issue of potential model misspecification. In particular, there is always the possibility of heterogeneous mean, beyond what is controlled by $x_i$ under the assumed model. As we discuss below, the matter affects the SP and QR approaches differently. Let $U_x = \{ i; x_i = x, i\in U\}$ be of the size $N_x$. Under the SP approach, which models $\mu_i$ by $\mu(x_i)$, heterogeneous mean is the case if $\mu_i \neq \mu(x_i)$, despite 
\begin{equation}
\mu(x) = \sum_{i\in U_x} \mu_i / N_x \label{h-mu}
\end{equation}
and the model is statistically correct in the sense that the $\mu_i$'s average to $\mu(x)$ for all the units in $U_x$. Notice that the condition \eqref{h-mu} can be verified in principle, provided non-informative $B$-sampling. Now that $E(y_i |\delta_i) = \mu_i$ by definition, we would have
\[
\sum_{i\in U_x} [E( y_i |\delta_i) - \mu(x)] = \sum_{i\in U_x} [\mu_i - \mu(x)] = 0
\]
Assuming $\mu_i = \mu(x)$ can still be prediction unbiased, despite heterogeneous mean. Meanwhile, under the QR approach, heterogeneous mean is the case if $p_i \neq p(x_i)$, despite 
\begin{equation}
p(x) = \sum_{i\in U_x} p_i / N_x  \label{h-p}
\end{equation}
Now that $E(\delta_i ; y_i) = p_i$ by definition, we would have
\[
E\big(\sum_{i\in U_x} \frac{\delta_i y_i}{p(x)} \big) - \sum_{i\in U_x}  y_i = p(x)^{-1} \sum_{i\in U_x} \big(p_i - p(x)\big) y_i \neq 0
\]
despite \eqref{h-p}. Insofar as such mean heterogeneity may be unavoidable in reality, the IPW estimator under the QP approach may be biased, even when the model of $p_i$ is statistically correct in the sense of \eqref{h-p}.

\section{More generally on validity conditions} \label{condition}

The discussion above suggests that the formulation of validity conditions in Section \ref{review} is inadequate in the presence of under-coverage and mean heterogeneity. Below we reformulate the validity conditions, which cover both the SP and QR approaches, despite the presence of under-coverage and mean heterogeneity. We elaborate and illustrate these conditions for the post-stratification and calibration estimators. Finally, we discuss the practical difficulties of verifying these validity conditions empirically.

\subsection{Non-parametric asymptotic (NPA) non-informativeness}

We start by noticing that, in the absence of any covariates, the $B$-sample mean equals to the population mean, denoted by $\bar{y}_B = \bar{Y}$, provided
\begin{gather*}
\begin{cases} Cov_N(\delta_i, y_i) = \frac{1}{N} \sum\limits_{i\in U} \delta_i y_i - 
\big( \frac{1}{N} \sum\limits_{i\in U} \delta_i \big) \big( \frac{1}{N} \sum\limits_{i\in U} y_i \big) = 0 \\
E_N(\delta_i) = \bar{p}_N \equiv \sum\limits_{i\in U} \delta_i/N > 0 \end{cases}
\end{gather*} 
where $E_N$ and $Cov_N$ denote, respectively, expectation and covariance with respect to the empirical distribution function that places point mass $1/N$ on each population unit. This is essentially an empirical formulation of the non-informativeness of the $B$-sample observation mechanism with respect to the outcome of interest. Similar expressions have appeared in various other discussions of the potential sample mean bias due to the observation mechanism, such as unequal probability sampling (Rao, 1966), survey nonresponse (Bethlehem, 1988), or big data (Meng, 2018). This motivates the following \textit{non-parametric asymptotic (NPA)} non-informativeness assumption in the absence of any covariates: 
\begin{equation}
\begin{cases}  \lim\limits_{N\rightarrow \infty} Cov_N(\delta_i, y_i) = 0 & \text{ i.e. non-informative B-selection} \\
 \lim\limits_{N\rightarrow \infty} E_N(\delta_i) = p >0 & \text{ i.e. non-negligible B-selection}\end{cases} \label{NPA}
\end{equation}
The NPA assumption \eqref{NPA} encompasses both the SP and QR approach. For the SP approach, taking the conditional expectation of $y_i$'s conditional on the $\delta_i$'s yields
\[
E\big(Cov_N(\delta_i, y_i) | \delta_U\big) = \frac{1}{N} \sum\limits_{i\in U} \delta_i \mu_i - 
\big( \frac{1}{N} \sum\limits_{i\in U} \delta_i \big) \big( \frac{1}{N} \sum\limits_{i\in U} \mu_i \big) \rightarrow 0
\]
provided NPA non-informative $B$-selection, where $\sum_{i\in U} \delta_i/N > 0$ given non-negligible $B$-selection in addition.
Under this condition, the $B$-sample expansion estimator \eqref{mean} is asymptotically prediction unbiased from the SP perspective. For the QR approach, taking the expectation of $\delta_i$'s with the $y_i$'s being constants yields
\[
\begin{cases} E\big(Cov_N(\delta_i, y_i) ; y_U\big) = \frac{1}{N} \sum\limits_{i\in U} p_i y_i - 
\big( \frac{1}{N} \sum\limits_{i\in U} p_i \big) \big( \frac{1}{N} \sum\limits_{i\in U} y_i \big) \rightarrow 0 \\
E\big( E_N(\delta_i) \big) = \sum_{i\in U} p_i/N \rightarrow p > 0 \end{cases}
\]
In particular, the NPA assumption \eqref{NPA} allows for $0\leq p_i\leq 1$, so that the $B$-sample expansion estimator \eqref{mean} remains consistent from the QR perspective, even in the presence of under-coverage of the units with $p_i =0$ or non-representative units with $p_i = 1$.

\subsection{Post-stratification estimator} 

Consider post-stratification by $x_i$, for $i\in U$. Provided the assumption \eqref{NPA} holds within each post-stratum, the $B$-sample post-stratification estimator is asymptotically unbiased from both the SP and QR perspective. Below we consider the matter under the QR approach. The SP approach is a special case of the calibration estimator discussed later.

Consider first the hypothetical  estimator with known $p_x = \sum_{i\in U_x} p_i/N_x$:
\[
\widetilde{Y} = \sum_x \sum_{i\in U_x} \delta_i y_i/p_x
\]
To fix the idea for variance estimation, suppose independent Bernoulli distribution of $\delta_i$ with probability $p_i$, where $0\leq p_i\leq 1$. The variance of $\widetilde{Y}$ is then given by
\[
V(\widetilde{Y}) =\sum_x \sum_{i\in U_x} p_i y_i^2/p_x^2  - \sum_x \sum_{i\in U_x} p_i^2 y_i^2/p_x^2 
\]
An unbiased estimator of the first term of the variance, denoted by $\tau_1$ is given by
\[
\hat{\tau}_1 = \sum_x \sum_{i\in U_x} \delta_i y_i^2/p_x^2 = \sum_x p_x^{-2} \sum_{i\in B_x} y_i^2
\]
where $B_x = B\cap U_x$. An unbiased estimator of the second term, denoted by $\tau_2$ is given by
\[
\hat{\tau}_2 = \sum_x p_x^{-2} \sum_{i\in U_x} \delta_i p_i y_i^2  
= \sum_x p_x^{-1} \sum_{i\in U_x} \delta_i y_i^2  = \sum_x p_x^{-1} \sum_{i\in B_x} y_i^2 
\]
where the second equality follows provide the additional QR$_x$ assumption, i.e. $p_i = p_x$ for $i\in U_x$. Putting $\hat{\tau}_1$ and $\hat{\tau}_2$ together, we obtain 
\[
\widehat{V}(\widetilde{Y}) = \sum_x \big( p_x^{-1} - 1 \big) p_x^{-1} \sum_{i\in B_x} y_i^2 
\]

Now, the post-stratification estimator, denoted by $\widehat{Y}$, is obtained from $\widetilde{Y}$ on replacing $p_x$ by $\hat{p}_x = n_{xB}/N_x$, where $n_{x,B}$ is the observed size of $B_x$. Its conditional variance given the observed $n_{xB}$'s for all $x$-values can be estimated by $\widehat{V}(\widetilde{Y})$ above, on replacing $p_x$ by $n_{xB}/N_x$. Provided $\widehat{Y}$ is asymptotically unbiased for $Y$, the conditional variance is approximately equal to the unconditional variance. Alternatively, expanding $\hat{p}_x$ around $p_x$ would yield an additional lower-order term due to $V(\hat{p}_x)$.

\subsection{Calibration estimator}

The post-stratification estimator is infeasible, in cases when the $B$-sample contains empty cells, or when the population size $N_x$ is not all known. Let 
\[
t_i = (t_{1i}, t_{2i}, ... t_{Ki})^{\top} = \big( t_1(x_i), t_2(x_i), ... t_K(x_i) \big)^{\top} = t(x_i)
\]
be a vector of many-to-one mappings of $x_i$, such that the population total $T = \sum_{i\in U} t_i$ is known, and the sample total $t = \sum_{i\in B} t_i$ has only non-zero components. 

As discussed for the calibration estimator in Section \ref{review}, generally one is not able to set the initial weight to be the inverse of $B$-sample inclusion probability in practice. Suppose one simply starts with initial equal weights $a_i = N/n_B$ for all $i\in B$. The linear calibration estimator (Deville and S\"{a}rndal, 1992) is given by
\[
\widehat{Y} = \sum_{i\in B} w_i y_i 
\]
where the weights $\{ w_i; i\in B\}$ minimise the distance to $\{ a_i; i\in B\}$ as measured by 
\[
\sum_{i\in B} (w_i - N/n_B)^2 
= \sum_t \big( \sum_{i\in B_t} w_i^2 - 2 (N/n_B) \sum_{i\in B_t} w_i + n_{tB} (N/n_{tB})^2 \big)
\]
subjected to the constraints $\sum_{i\in B} w_i t_i = T$, where $B_t = \{ i; t_i = t, i\in B\}$ and $n_{tB} >0$. It follows that $w_i = w_t$, for $i\in B_t$, since (i) the only thing that matters to the calibration constraints is $\sum_{i\in B_t} w_i$ now that $t_i = t$ for $i\in B_t$, and (ii) given whatever $\sum_{i\in B_t} w_i$, the term $\sum_{i\in B_t} w_i^2$ is minimised at $w_i = w_t$ for $i\in B_t$. 

As the first validity condition for $\widehat{Y}$, suppose there exists a vector $\beta_{K\times 1}$, such that 
\begin{equation}
\sum_{i\in U_t} \epsilon_i/N_t \rightarrow 0 \label{NPA-SP}
\end{equation}
for any possible $t$, as $N\rightarrow \infty$, where $\epsilon_i = y_i - t_i^{\top} \beta$, and $N_t$ is the population size of $U_t = \{ i; t_i = t, i\in U\}$. We shall refer to \eqref{NPA-SP} as the NPA$_t$ assumption, which is essentially an NPA version of the SP$_x$ assumption \eqref{SPx}, where the covariate $x_i$ is replaced by $t_i$ here. Moreover, it relaxes the model \eqref{SPx} of the conditional mean, allowing for potential heterogeneous mean similar to \eqref{h-mu}. Now that $\sum_{i\in B} w_i t_i = T$, we have 
\[
\widehat{Y} - Y = \sum_{i\in B} w_i (t_i^{\top} \beta + \epsilon_i) - \sum_{i\in U} t_i^{\top} (\beta + \epsilon_i)
= \sum_{i\in B} w_i \epsilon_i - \sum_{i\in U} \epsilon_i \\
\]
Provided the NPA$_t$ assumption, $\sum_{i\in U} \epsilon_i/N \rightarrow 0$ as $N\rightarrow \infty$. Moreover, we have
\[
\frac{1}{N} \sum_{i\in B} w_i \epsilon_i = \sum_t \frac{w_t}{N} \sum_{i\in U_t} \delta_i \epsilon_i 
= \sum_t w_t \frac{N_t}{N} \big( Cov_{N_t}(\delta_i, \epsilon_i) + (\frac{1}{N_t} \sum_{i\in U_t} \delta_i) (\frac{1}{N_t} \sum_{i\in U_t} \epsilon_i) \big)
\rightarrow 0
\]
as $N\rightarrow \infty$, provided 
\begin{equation}
\begin{cases} Cov_{N_t}(\delta_i, \epsilon_i) \rightarrow 0 \\ E_{N_t}(\delta_i) = \sum_{i\in U_t} \delta_i/N_t \rightarrow p_t > 0 \end{cases} \label{NPA-niy}
\end{equation} 
for any given $t$, which is an adaption of the NPA non-informativeness assumption \eqref{NPA} to the present setting. It follows that the two NPA assumptions \eqref{NPA-SP} and \eqref{NPA-niy} are the key validity conditions for the calibration estimator to be consistent for $Y$.

For variance estimation, suppose again independent Bernoulli distribution of $\delta_i$ with probability $p_i$, where $0\leq p_i\leq 1$. An approximate variance estimator for the calibration estimator $\widehat{Y}$ can then be given as
\[
\widehat{V}(\widehat{Y}) = \sum_t \big( \hat{p}_t^{-1} - 1 \big) \hat{p}_t^{-1} \sum_{i\in B_t} (y_i - t_i^{\top} \hat{\beta})^2 
\]
where $\hat{p}_t = n_{tB}/N_t$, and $\hat{\beta} = \big( \sum_{i\in B} w_i t_i t_i^{\top} \big)^{-1} \big( \sum_{i\in B} w_i t_i y_i \big)$.

\subsection{Validation of non-informative $B$-sample selection}
 
Of the validity conditions discussed above, the critical assumption is non-informative $B$-sample selection, which can be stated in various forms. For instance, provided the non-informativeness assumption \eqref{NPA}, an additional assumption like \eqref{NPA-SP} can in principle to validated empirically. However, the non-informativeness condition may not hold exactly, and it is generally impossible to verify only based on the data used for the estimation. Below we discuss the issue in more details.

Consider first the setting without any relevant covariates. Let $\hat{p}_i$ be any possible estimator of $p_i = E(\delta_i ; y_i)$, for $i\in U$. One has only two empirical checks for them:
\[
\begin{cases} \sum\limits_{i\in B} 1/\hat{p}_i = N\\  \sum\limits_{i\in U} \hat{p}_i = n_B \end{cases} 
\]
Both of them are perfectly satisfied by $\hat{p}_i \equiv n_B/N$. Since assuming constant $p_i$ would satisfy the NPA non-informative assumption \eqref{NPA} trivially, it is not possible to use these checks to diagnostic potential departures from the assumption \eqref{NPA}. Clearly, the same difficulty exist for the post-stratification estimator.

Consider next the propensity model $p_i = p(x_i; \lambda)$, with covariates $x_i$ and parameter $\lambda$.
Provided known $x_U$, the census score equation is given by
\[
\sum_x \frac{\partial p(x;\lambda)}{\partial \lambda} \left[ \frac{n_{xB}}{p(x;\lambda)} - \frac{N_x - n_{xB}}{1- p(x;\lambda)} \right] =0  
\] 
which is always satisfied by $p(x;\hat{\lambda}) = n_{xB}/N_x$, i.e. the saturated model. It follows that for any non-saturated model, the potential lack-of-fit can always be attributed to the misspecification of the functional form $p(x_i; \lambda)$, but not that $p_i$ can be given as a function of $x_i$. Thus, the validity of propensity modelling cannot be refuted empirically. 

Finally, assume the $B$-sample inclusion probability $p_i$ depend on $x_i$, where $x_i$ is known for $i\in U$. For goodness-of-fit checks, let $z_i$ be a known covariate, which is distinct from $x_i$ and cannot be obtained from $x_i$ via a many-to-one mapping. We have
\[
\begin{cases} E(z_B) = \sum_{i\in U} p_i z_i = \sum_x p(x; \lambda) \sum_{i\in U_x} z_i = \sum_x p(x;\lambda) N_z \bar{Z}_x\\
Z = E(\sum_{i\in U} \delta_i z_i/p_i) = E[\sum_x n_{xB} \bar{z}_{xB}/p(x;\lambda)] \end{cases}
\]
where $\bar{Z}_x = \sum_{i\in U_x} z_i/N_x$ and $\bar{z}_{xB} = \sum_{i\in B_x} z_i/n_{xB}$. The two observed checks are 
\[
\begin{cases} z_B \equiv \sum_x n_{xB} \bar{z}_{xB} = \sum_x \hat{p}_x N_x \bar{Z}_x \\ Z = \sum_x n_{xB} \bar{z}_{xB}/\hat{p}_x \end{cases}
\]
Setting $\hat{p}_x = n_{xB}/N_x$, which fits the assumption $p_i = p(x_i; \lambda)$, the two checks are satisfied provided $\bar{Z}_x = \bar{z}_{xB}$, i.e. the $B$-sample expansion estimate of $Z_x$ is perfect for all $x$. This would suggest that the NPA assumption \eqref{NPA} holds for $z_i$ given $x_i$. This may be considered to support the plausibility of the NPA assumption \eqref{NPA} for $y_i$ given $x_i$, provided $z_i$ is known to be correlated with $y_i$, but not otherwise. However, in situations where such a covariate $z_i$ is available, it seems natural that it should be used in the estimation of $Y$ to start with. Thus, one is faced with a dilemma, where building the best model for estimation would at the same time reduce the ability to verify it.

\section{Using additional probability sample of outcomes} \label{combine} 
  
So far we have only considered the situations, where the outcome values of interest are only observed in the non-probability sample $B$. As discussed in Section \ref{review}, the availability of relevant covariates in a second probability sample $S$ cannot help to address the conceptual validity of inference, although it may help to overcome the practical difficulty when these covariates are not known for the whole population. Obviously, the matter changes completely, provided instead a second probability sample of outcomes. Below we discuss shortly two different approaches to inference in the absence of any relevant covariates. The ideas remain the same in situations with additional covariates. 

To start with, given the non-probability sample observations $y_B$, one may treat $(B, y_B)$ as fixed, and select a second supplementary sample from the rest of the population, denoted by $S\subset U\setminus B$. Provided the $S$-sample observations of the outcome, denoted by $y_S$, it is straightforward to obtain a test for $H_0: \bar{Y} = \bar{y}_B$, given as
\[
\eta = (\bar{y}_B - \widehat{\bar{Y}_B^c})^2 / \widehat{V}(\widehat{\bar{Y}_B^c})
\] 
where $\widehat{\bar{Y}_B^c}$ is an $S$-sample estimator of the population mean outside of the $B$-sample, i.e.
\[
\bar{Y}_B^c = \sum_{i\in U\setminus B} y_i/(N - n_B) 
\]
and $\widehat{V}(\widehat{\bar{Y}_B^c})$ is the associated variance estimator. In many situations, one may need to reject $H_0$, provided small enough $\widehat{V}(\widehat{\bar{Y}_B^c})$. Provided this is the case, let $W_B = n_B/N$. A consistent estimator of $\bar{Y}$ is then given by
\[
\widehat{\bar{Y}}_S = W_B \bar{y}_B + (1 - W_B) \bar{y}_w \qquad\text{and}\qquad 
\bar{y}_w = \frac{\sum_{i\in S} y_i/\pi_i}{\sum_{i\in S} 1/\pi_i}
\]
where $\pi_i$ is the $S$-sample inclusion probability, and the validity of $\widehat{\bar{Y}}_S$ now derives from probability sampling of $S$, regardless of how the $B$-sample is generated. The relative efficiency (RE) against the setting without the $B$-sample can be given by
\[
\mbox{RE} = \big[ (1-W_B)^2 V(\bar{y}_w) \big]/V(\widehat{\bar{Y}'}) 
\]
where $\widehat{\bar{Y}'}$ is a hypothetical probability sample from the whole population $U$, which has the same sample size and the same sampling design as $S$. One may refer to this as the \emph{split-population} approach to inference, which is an age-old idea for combining survey sampling with administrative data. The efficiency gain would be substantial provided the $B$-sample is large. In fact, the larger the $B$-sample, the greater is the efficiency gain.
 
Under the second approach to inference, consider a \emph{composite} estimator given by  
\[
\widehat{\bar{Y}}_C = \gamma \bar{y}_B + (1 - \gamma) \bar{y}_w 
\]
where $\gamma$ is the composition weight, for $W_B \leq \gamma \leq 1$. Notice that when $\gamma = W_B$, the composite estimator is just the split-population estimator $\widehat{\bar{Y}}_S$ above, which is consistent for $\bar{Y}$. As $\gamma$ increases from $W_B$ towards one, one risks introducing greater bias, insofar as $\bar{y}_B \neq \bar{Y}$. However, the composite estimator may yield a smaller mean squared error (MSE) of estimation, provided this is desirable. One is then essentially trading the increasing bias $(\gamma - W_B) (\bar{y}_B - \bar{Y}_B^c)$ against the decreasing stand error $(1- \gamma) \mbox{SE}(\bar{y}_w)$, as $\gamma$ increases.
The composite estimator that achieves the minimum MSE is given by 
\[
\gamma = \frac{V(\bar{y}_w) + W_B (\bar{y}_B - \bar{Y}_B^c)^2}{V(\bar{y}_w) + (\bar{y}_B - \bar{Y}_B^c)^2} 
\]
Estimating $\bar{Y}_B^c$ by $\bar{y}_w$ in application, one can use
\[
\hat{\gamma} = \min(W_B + (1-W_B) \widehat{V}(\bar{y}_w)/(\bar{y}_B - \bar{y}_w)^2, 1)
\]
The validity of the composite approach derives from probability sampling of $S$, regardless of how the $B$-sample is generated. Again, the bigger the $B$-sample, the better it is.

\end{document}